\documentclass[12pt]{article}

\usepackage[tbtags]{amsmath}
\usepackage{amsfonts,amssymb}

\usepackage{mathrsfs}
\usepackage{graphicx}
\usepackage{epstopdf}

\usepackage{multirow}

\newtheorem{statement}{Statement}

\newtheorem{proof}{Proof}

\providecommand{\keywords}[1]{\textbf{\textit{Key words:}} #1}

\begin{document}
\author{Vladimir\,L.~Gavrikov}
\title{Some features of natural sequence totals in various numeral systems}
\maketitle

\begin{abstract}
Minor totals of natural sequence were shown to possess some properties in respect to their units digits. Depending on numeral system applied the units digits may take any digit of the system or there may be exclusions. i.e. some system digits can never appear in the units digits of the minor totals. Natural sequence arithmetics may be applied to model multicellular organism growth. It can be shown that growth of an idealized organism may follow a power law in the form of the minor totals.
\end{abstract}
\\

Let's consider first-order arithmetical progression such as natural sequence in the form $1 + 2 + 3 + 4 + \dots$. The formula of minor totals of the natural sequence has been known since long ego (there are evidences that Pythagoras had known about it \cite{Pengelley}). The $n_{th}$ minor total of the sequence $S_n$ is

\begin{equation}
S_n = \sum^n_{\kappa =1}\kappa = \frac{n(n+1)}{2},
\label{eq:1}
\end{equation}
and any number of minor totals may be calculated with the formula.

Basic observations over $S_n$ give an idea that the set of these minor totals show some patterns. For instance, in base-ten system the units digits of $S_n$ ($n$ being large enough) are found with uneven frequencies. Namely, $0 (0.2)$, $1 (0.2)$, $2 (0)$, $3 (0.1)$, $4 (0),$ $5 (0.2)$, $6 (0.2)$, $7 (0)$, $8 (0.1)$, $9 (0)$, where $0, 1, 2, 3, 4, 5, 6, 7, 8, 9$ are digits of base-ten system and approximate empirical frequencies of appearance of the units digits of $S_n$ are given in brackets. In other words, the minor totals $S_n$ under any $n$ end never (empirically) with digits $2, 4, 7, 9$. Let's denote it in such a way that the frequency distribution is gappy, i.e. contains gaps which are zeros of some digit's frequencies.

The other result one gets considering $S_n$ in base-eight system. Empirically, all the digits of the base-eight system $0$, $1$, $2$, $3$, $4$, $5$, $6$, $7$ appear with approximately the same frequencies in order of units of minor totals $S_n$, i.e. the frequency distributions of the digits are continuous (non-gappy).

Considering then numbering systems with other bases one can find that gappy frequency distributions of units digits of $S_n$ appear in systems with $3$, $5$, $6$, $7$, $9$, $10$ bases, while non-gappy distributions are found in system with $4$, $8$, $16$ bases. In other words, one could imagine that the systems with $4$, $8$, $16$ bases are different from other system as long as the matter is the content of units digits of $S_n$.

The above considerations may be illustrated on the example of base-three and base-four systems.

\begin{statement}
In base-three system, the minor totals of the natural sequence $S_n$ never end with $2$ in the order of units.
\end{statement}

\begin{proof}
Consider minor totals on the natural sequence in the form of equations
\begin{equation}
S_{3k+i} = 3m + j,
\label{eq:2}
\end{equation}
where expressions $3k + i$ and $3m + j$ show the units digits ($i$ and $j$) in base-three representation of the length of natural sequence ($n$) and its minor total ($S_n$), corres\-pondingly. Obviously, $i, j \in \{ 0, 1, 2 \}$ in base-three system. By definition, $k$ and $m$ are whole numbers, which is of key importance in the proof.

All the variety of the natural sequence lengths presented as $3k + i$ at any $k$ is exhausted by $3k + 0$, $3k + 1$, $3k + 2$. Let's consider them consecutively. To be more precise, let's check which $j$ (three-base units digits) are possible in the equations (\ref{eq:2}).

First, let $n = 3k + 0$. Then, according to (\ref{eq:1}), $S_{3k+0} = \cfrac{3k(3k+1)}{2}$. Let's check which of equations $\cfrac{3k(3k+1)}{2} = 3m + 0$, $\cfrac{3k(3k+1)}{2} = 3m + 1$, $\cfrac{3k(3k+1)}{2} = 3m + 2$ are correct in that sense that they do not violate the condition ``$k$ and $m$ are whole numbers'':
\begin{eqnarray}
j = 0 & \Rightarrow & m = \frac{k(3k + 1)}{2}, \nonumber \\
j = 1 & \Rightarrow & m = \frac{k(3k + 1)}{2} - \frac{1}{3}, \label{eq:3} \\
j = 2 & \Rightarrow & m = \frac{k(3k + 1)}{2} - \frac{2}{3}. \nonumber
\end{eqnarray}

Because at any whole $k$ the expression $k(3k + 1)$ is always even then $\cfrac{k(3k + 1)}{2}$ is always a whole number. It is obvious thus that only $j = 0$ in (\ref{eq:3}) does not violate the initial conditions. In other cases $m$ appears to be a fractional number; a fractional $m$ leads to that numbers appear among $S_n$ that do not belong to minor totals in reality. For example, at $j = 1$ and $k = 5$ in (\ref{eq:3}) $S_n$ would equal 119 but this number does not belong to the set of the minor totals of $S_n$ according to (\ref{eq:1}). Thus at $n = 3k + 0$ the correspondent minor total $S_n$ in three-base representation can have only $0$ in units digit.

Then let $n = 3k + 1$. Analogically, it should be checked which $j$ do not lead to contradictions. So, $S_{3k+1} = \cfrac{(3k+1)(3k+2)}{2} = 3m + j$ and:

\begin{eqnarray}
j = 0 & \Rightarrow & m = \frac{3k(k + 1)}{2} + \frac{1}{3}, \nonumber \\
j = 1 & \Rightarrow & m = \frac{3k(k + 1)}{2}, \label{eq:4} \\
j = 2 & \Rightarrow & m = \frac{3k(k + 1)}{2} - \frac{1}{3}. \nonumber
\end{eqnarray}

At any whole $k$, the expression $\cfrac{3k(k + 1)}{2}$ is always a whole number. From (\ref{eq:4}), it is seen that the only case of $j = 1$ does not lead to contradictions. Thus if the natural sequence length is $n = 3k + 1$ then its minor total in three-base expression has $1$ in units digit.

Lastly, let $n = 3k + 2$. In this case $S_{3k+2} = \cfrac{(3k+2)(3k+3)}{2} = 3m + j$ and:

\begin{eqnarray}
j = 0 & \Rightarrow & m = \frac{k(3k + 5)}{2} + 1, \nonumber \\
j = 1 & \Rightarrow & m = \frac{k(3k + 5)}{2} + \frac{2}{3}, \label{eq:5} \\
j = 2 & \Rightarrow & m = \frac{k(3k + 5)}{2} + \frac{1}{3}. \nonumber
\end{eqnarray}

It can be seen from (\ref{eq:5}) that the only case $j = 0$ does not lead to contradictions because $\cfrac{k(3k + 5)}{2}$ is alway a whole number at any whole $k$.

Therefore all the possible combinations of $i$ and $j$ have been considered in (\ref{eq:3}) -- (\ref{eq:5}) and none of these combinations leads to that three-based minor total $S_n$ has $2$ in units digit.
\end{proof}

\begin{statement}
In four-base numbering system, the minor totals of the natural sequence $S_n$ can have any digit of the system ($0$ or $1$ or $2$ or $3$) in their units digits.
\end{statement}

\begin{proof}
Analogically to the case given above, let's consider minor totals of natural sequence in the form of the equations
\begin{equation}
S_{4k+i} = 4m + j,
\label{eq:6}
\end{equation}
where $4k + i$ and $4m + j$ are expressions showing units digits $i$ and $j$ in the length of natural sequence $n$ and its minor total $S_n$, correspondingly, in four-base system.

Let' s consider the combinations between $i$ and $j$ sequentially.

First, $n = 4k + 0$. Then from (\ref{eq:1}) the minor total $S_{4k+0} = \cfrac{4k(4k + 1)}{2}$ and considering all the $j \in \{ 0, 1, 2, 3 \}$ one can get:
\begin{eqnarray}
j = 0 & \Rightarrow & m = \frac{k(4k + 1)}{2}, \label{eq:7.1} \\
j = 1 & \Rightarrow & m = \frac{2k(4k + 1) - 1}{4} \label{eq:7.2}, \\
j = 2 & \Rightarrow & m = \frac{k(4k + 1) - 1}{2} \label{eq:7.3}, \\
j = 3 & \Rightarrow & m = \frac{2(4k + 1) - 3}{4} \label{eq:7.4}.
\end{eqnarray}

The expression for $m$ in (\ref{eq:7.1}) is not correct for odd $k$ because it brings about a fractional $m$. Still for even $k$ this expression is correct, i.e. the minor total in four-base system can have $0$ in units digit.

Because $2k(4k + 1) - 1$ in (\ref{eq:7.2}) is always odd then $m$ is always fractional. The same is true for $2k(4k + 1) - 3$ in (\ref{eq:7.4}). In (\ref{eq:7.3}), the expression $k(4k + 1) - 1$ is even at odd $k$ and thus $m$ is  whole in this case.

Consequently, at $n = 4k + 0$ the four-base minor total of natural sequence can have $0$ and $2$ in units digit.

Second, let $n = 4k + 1$. From (\ref{eq:1}) the minor total\\ $S_{4k+1} = \cfrac{(4k + 1)(4k + 2)}{2}$. Enumera\-tion of possibilities for $j \in \{ 0, 1, 2, 3 \}$ looks like:
\begin{eqnarray}
j = 0 & \Rightarrow & m = \frac{k(4k + 3)}{2} + \frac{1}{4}, \label{eq:11} \\
j = 1 & \Rightarrow & m = \frac{k(4k + 3)}{2} \label{eq:12}, \\
j = 2 & \Rightarrow & m = \frac{k(4k + 3)}{2} - \frac{1}{4}, \label{eq:13} \\
j = 3 & \Rightarrow & m = \frac{k(4k + 3)}{2} - \frac{1}{2}. \label{eq:14}
\end{eqnarray}

The possibilities shown in (\ref{eq:11}) and (\ref{eq:13}) lead to that $m$ is fractional at any whole $k$. In (\ref{eq:12}), $m$ is fractional at odd $k$ while $m$ is whole at even $k$. In (\ref{eq:14}), $m$ is fractional at even $k$ while $m$ is whole at odd $k$. Thus in the case of $n = 4k + 1$ the four-base minor totals of natural sequence can have $1$ and $3$ in units digits.

Therefore while considering only a part of combinations between $i$ and $j$ one can see that the four-base minor totals can have all the $\{ 0, 1, 2, 3 \}$ in units digits.
\end{proof}

The statements above may be reformulated in terms of modular arithmetics. So, $S_{3k+0} \equiv 0$ (mod $3$), $S_{3k+1} \equiv 1$ (mod $3$), $S_{3k+2} \equiv 0$ (mod $3$) while $S_n \not\equiv 2$ (mod $3$) at any $n$. Then, $S_{4k+0} \equiv 0$ (mod $4$) at even $k$, $S_{4k+0} \equiv 2$ (mod $4$) at odd $k$, $S_{4k+1} \equiv 1$ (mod $4$) at even $k$, $S_{4k+1} \equiv 3$ (mod $4$) at odd $k$.

Through applying of the approach given above one can prove presence or absence of gaps in frequency distributions of units digits of minor totals $S_n$ of natural sequence under numbering system with other bases. Particularly, the gaps are found in numbering systems with bases $5$, $6$, $7$, $9$, $10$ while there are no gaps in systems with bases $8$, $16$. For the proving, it is necessary to enumerate all the combinations between $i$ and $j$ in, correspondingly, lengths of natural sequence $n = Lk + i$ and its minor totals $S_n = Lm + j$, $L$ being the base of numbering system.

\begin{figure}[tb]
\center{\includegraphics[width=0.6\textwidth]{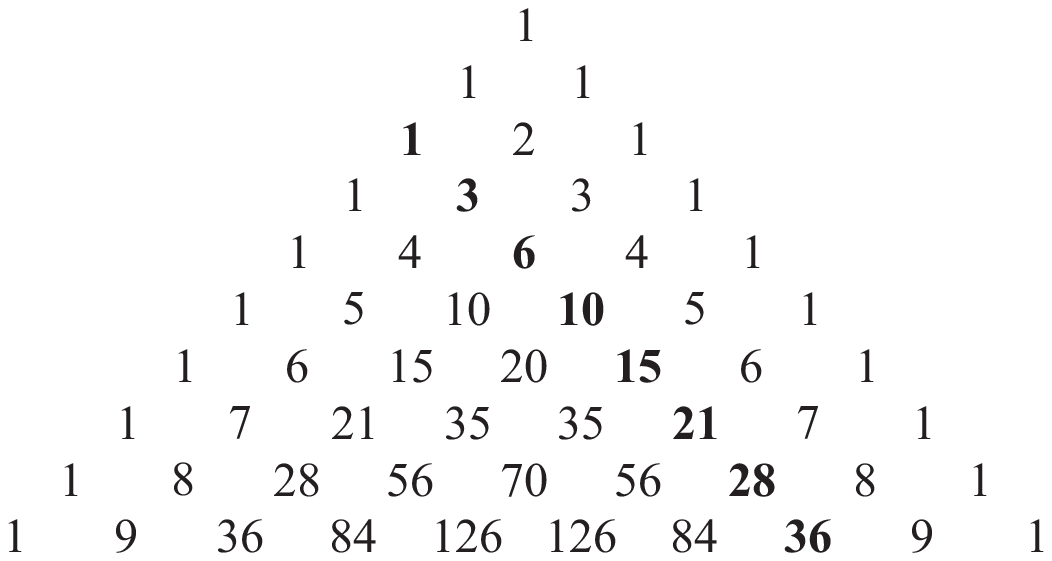}}
\caption{Pascal triangle. Minor totals of natural sequence are given in bold face.}
\label{fig1}
\end{figure}

The minor totals (\ref{eq:1}) considered above have a relation not only to the natural sequence as such. As it is well known, natural sequence and its minor totals are particular cases in the system of binomial coefficients. To receive evidence of that it is enough to have a look at Pascal triangle (fig. \ref{fig1}) which is a triangular table where every $n_{th}$ row is coefficients of binomial theorem $(1 + x)^n$ ($x$ being a real variable). It is easy to see that the minor totals of natural sequence are coefficients at the quadratic terms ($n\geq 2$).

Arithmetical properties of binomial coefficients are known since long ego and relate mostly to the divisibility by prime numbers and their degrees \cite{Winberg} (Lucas, Kummer theorems etc.) as well as to sums of the coefficients in one row. As it can be seen from the analysis given above some subsets of the coefficients from diagonals of the Pascal triangle can possess other properties as well. Particularly, the frequency distributions of units digits of binomial coefficients at quadratic terms (which are minor totals of natural sequence) may or may not have gaps dependently on numbering system considered. This property can be strictly proved.

\begin{figure}[tb]
\center{\includegraphics[width=0.6\textwidth]{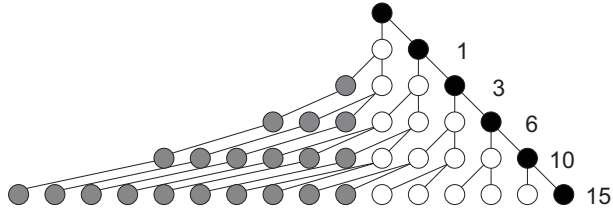}}
\caption{A schematic idealization of cell population growth within a multicellular organism. Black filled circles denote initial cells, open circles are dividing descendant cells, gray filled circles are differentiated (non-dividing) cells. Numbers to the right show quantities of all the descendant cells.}
\label{fig2}
\end{figure}

The arithmetics of natural sequence has also applications in the area of biological combinatorics or, to be more precise, in modeling of biological growth. Size of an organism is a function of number of cells in it and their sizes. While cell sizes are continuous quantities and require correspondent mathematics the numbers of cells are whole numbers. Thus in terms of cell number an organism growth may be considered with the help of whole number mathematics.

A classical example of whole number modeling is provided by a consideration of population growth of cells each of which doubles in equal time intervals producing the next cell generation. Such a dynamics gives rise to sequence of the cell numbers like $1$, $2$, $4$, $8$, $\dots$ , $2^i$, $i$ being the number of the cell generation while $2$ in this particular case being the reproduction factor. In other words, this is an example how local divisions of cell lead to an \textit{exponential law} of cell population growth.

Within a multicellular organism, however, the cell population growth obeys much more complicated rules. On the example of higher plant organisms the following rules (kinds of cells) may be identified. First, there are the so-called initial cells that preserve the ability to divide in the course of the entire life span of the organism. Second, the cells--immediate descendants of the initials--can divide for some time but sooner or later transform to the third kind of cells. The third kind are remote descendants that are differentiated cells that may be dead or alive but no longer divide.

Suppose, the growth goes like that depicted in fig. \ref{fig2}, i.e. the initial cell reproduces itself and gives life to a descendant cell. The descendant cell produces other descendant cells that can produce both descendant cells and differentiated cells. Suppose next that the number of dividing cells grows linearly constituting a sort of natural sequence, which is quite plausible for the initial stages of growth.

Calculating the total quantities of descendant cells (immediate plus differentiated descendants) one can easily see (fig. \ref{fig2}) that the cell population grows as minor totals of the natural sequence. A simple corollary of the dynamics is that the growth of this idealized multicellular organism obey a \textit{power law} which in the considered pure case is given by (\ref{eq:1}). It is easy to show that if the number of dividing cells stays constant then the entire organism grows linearly, in terms of cell quantity. If the number of dividing cells falls linearly then the entire organism undergoes a decay of growth which obeys a power law as well. The question of why the number of dividing cells follow this or that dynamics is not considered here. The important implication of the cell population behavior modeled with natural sequence arithmetics is that it is a power law that govern growth of a multicellular organism, naturally, in terms of cell number.

Another implication of cell division modeling is that if the process follows the logic of natural sequence summation then the number of cell might obey the same properties of minor totals that were considered above. Particularly, depending on the numbering system applied the number of cell may possess only certain values, which leads to predictable gaps in frequency distributions of cell numbers in the course of growth.\\

Siberian federal university, Krasnoyarsk, Russia

vgavrikov@sfu-kras.ru

\date{}

\end{document}